\documentstyle{amsppt}

\input amstex
\NoBlackBoxes 
\topmatter
\title
Quasideterminants, I
\endtitle
\author
I. Gelfand and V. Retakh
\endauthor

\address
\newline
I.~G.: Department of Mathematics, Rutgers University,
New Brunswick, NJ 08903
\newline
V.~R.: Department of Mathematics, Harvard University,
Cambridge, MA 02138
\endaddress
\email
\newline
I.~G.: igelfand\@math.rutgers.edu
\newline
V.~R.: retakh\@math.harvard.edu
\endemail
\endtopmatter
\centerline{\bf Introduction}

A notion of quasideterminants for matrices over a noncommutative
skew-field was introduced in [GR], [GR1], [GR2]. It proved its
effectiveness for many areas including noncommutative symmetric
functions [GKLLRT], noncommutative integrable systems [RS],[EGR],
quantum algebras [GR], [GR1], [GR2], [KL], [M], etc.

The main property of quasideterminants is a "heredity principle":
let $A$ be a square matrix over a skew-field and $(A_{ij})$ be its
block decomposition into submatrices of $A$. Consider $A_{ij}$'s 
as elements of a matrix $X$. Then a quasideterminant of matrix X will
be a matrix  $B$ again, and (under a natural assumption) a 
quasideterminant of $B$ will be equal to a suitable quasideterminant
of $A$.

This principle is not valid for commutative determinants, 
because they are not defined for block-matrices.

Quasideterminants are not analogues of commutative determinants but
rather a ratio of determinants of $n\times n$-matrices to
determinants of their $(n-1)\times (n-1)$-submatrices. In fact,
for a matrix over a commutative algebra a quasideterminant is 
equal to a correspondence ratio.

Many noncommutative areas of mathematics (Ore rings, rings of 
differential operators, theory of factors, "quantum mathematics",
Clifford algebras, etc) were developed separately from each other.
 Our approach shows an advantage of working with totally 
noncommutative variables (over free rings and skew-fields).
 It leads us to a big variety
of results, and their specialization to different noncommutative
areas implies known theorems with additional information.

The price you pay for this is a huge number of inversions in
rational noncommutative expressions measured by their height which cannot
be reduced. But this invariant (inversion height) shows us a "degree of 
noncommutativity" and it is of a great interest by itself.

Our experience shows that dealing with noncommutative objects one
should not imitate the classical commutative mathematics, but follow
"the way it is" starting with basics. In this paper we
consider mainly two such problems: noncommutative Pl\"ucker coordinates
(as a background of a noncommutative geometry) and  
noncommutative Bezout and Vieta theorems (as a background of noncommutative
algebra). We apply these results to the theory of noncommutative symmetric
functions started in [GKLLRT].

The first chapter of this paper contains basic
definitions and properties of quasideterminants. We consider
Pl\"ucker coordinates in Chapter II. Chapter III is devoted to
Bezout and Vieta theorems and also noncommutative symmetric
functions. In the last Chapter we continue our investigation
of noncommutative continued fractions and almost triangular matrices.
 It turns out that this problem
is related with computation of quantum cohomology [G], [FK].  

Some results in the paper are given without proofs, since proofs
are straightforward or contained in previous publications.

We thank P. Etingof for his careful reading of the manucript.

\head  I.  A general theory and main identities\endhead
\medskip
{\centerline {\bf \S 1.1 Definition}}

\subhead  I.1.1.  The skew-field of rational functions of free
variables\endsubhead
\medskip
This theory was developed by Amitsur, Bergman and P. M. Cohn (see, for
example [C], [C1]). We consider here skew-fields over fields of
characteristic $0$, but everything may be generalized for fields
of characteristics $p$. 
Here we will remind the approach by Amitsur presented by 
Bergman.  Given a set $X =
\{x_1,\ldots,x_n\}$ we write $\Bbb F(X)$ for a the algebra freely
 generated from the elements of $X$ over rationals $\Bbb Q$ by
the operations of addition, subtraction, multiplication and taking an
inverse.
\medskip
No relations are imposed, thus e.g. $(x - x)^{-1}$ exists.  The
algebra $\Bbb F(X)$ is called an algebra of rational formulas on $X$.
Let $R$ be a ring with unit.  Any map $\alpha: X\rightarrow R$ defines
a map $\overline\alpha$ of a subset of $\Bbb F(X)$ into $R$ by the
following rules
\roster
\item"{i)}"  $\overline\alpha (m) = m, \ m\in \Bbb Z$,

\item"{ii)}"  $\overline\alpha (x_i) = \alpha(x_i), \quad i
= 1,\ldots,n$,

\item"{iii)}"  if $a =-b$, or $b + c$, or $bc$ and
$\overline\alpha(b), \overline\alpha(c)$ are defined then
$\overline\alpha(a) = - \overline\alpha(b)$, or $\overline\alpha (b +
c)$, or $\overline\alpha(b) \cdot\overline\alpha (c)$,

\item"{iv)}"  if $a = b^{-1}$ and $\overline\alpha(b)$ is defined and
invertible in $R$ then $\overline\alpha(a) = (\overline\alpha
(b))^{-1}$.
\endroster
\medskip
Let now $\alpha : X\rightarrow D$ be a map to a skew-field $D$.  Then
$\alpha(a)$ is undefined if and only if $a$ has a subexpression
$b^{-1}$ such that $\alpha(b) = 0$.  For each $\alpha : X\rightarrow
D$ one can consider a subset $E(\alpha )$ of $\Bbb F(X)$, the domain of
$\overline\alpha$, consisting of the expressions which can be
evaluated for $\alpha$.  Similarly with each $f\in\Bbb F(X)$ one can 
associate its domain dom$f$, a subset of $D^n$ consisting of the
points $(\alpha =(\alpha _i), i=1,\dots n)$ such that$f\in E(\alpha )$ 
is defined; $f$ is called {\it nondegenerate} if
dom$f\not=\emptyset$.
\medskip
One can show [C],[C1] that if $D$ is a skew-field which is an algebra
over an infinite field $k$ and $f$ and $g$ are nondegenerate then
dom$f\cap $dom$g\not=\emptyset$.
\medskip
Given $f,g\in\Bbb F(X)$, let us put $f\sim g$ if $f,g$ are
nondegenerate and $f,  g$ have the same value at each point of
dom $f\cap $dom$g$.  This is an equivalence.
\medskip
We call the classes of equivalence the rational functions on a
skew-field $D$.  Then we have
\medskip 
\proclaim{Theorem 1.1.1 [C], [C1]} 

$i)$ Let $D$ be a skew-field with a
center $k$ of characteritic $0$.  Then the equivalence classes of 
rational formulas
of variables $x_1,\ldots,x_n$ form a skew-field $F_D(X)$.

$ii)$ If $D$ is infinite dimensional over $k$ then $F_D(X)$ does not
depend on the skew-field $D$.
\endproclaim
Under assumptions of i) and ii) we will identify all $F_D(X)$ and
use a notation $F(X)$. For example, if $\# X=1$, i.e. $X={x}$ then
$F(X)=\Bbb Q(x)$.  

This skew-field $F(X)$ is called {\it the free skew-field}.

All the equalities in this chapter are considered as the equalities of
rational functions. 
\remark{Remark} As shown by P.M. Cohn [C], [C1] the free field is
universal in the following sense: consider the ring $k\langle
x=(x_1,\dots ,x_n)\rangle$ of noncommutative polynomials over the 
commutative field
$k$, and the category whose objects are homomorphisms $k\langle
x\rangle \to P, P$ is a skew-field; the morphisms between two such
objects $\phi _i: k\langle x \rangle \to P_i, i=1,2$ are 
{\it specializations} $s: P_1\to P_2$, i.e. homomorphisms from
a local subring of $P_1$ into $P_2$ such that $s\phi _1=\phi _2$
(see [C], [C1] for exact definitions).
Then the free skew-field is an initial object in this category.
\endremark
We remind the following definition (see, for example, [R]).
\definition{Definition 1.1.2}  The inversion height of a rational 
expression is the
maximal number of nested inversions in it.  The inversion height of an
element in the free skew-field $F(X)$ is the smallest inversion height of its
rational expressions.
\enddefinition

Let $D$ be a skew-field generated by elements $x_1,\dots ,x_n$ over
a central subfield. If $y\in D$ may be written as a polynomial
expression $P(x_1,\dots ,x_n)$ in $D$ then $y$ has inversion height
zero. If $x_1,\dots ,x_n$ are generators of an Ore ring, and $D$
its skew-field of fractions then any element $z\in D$ is a ratio of
two polynomial expressions of generators. Then, $z$ has height
$\leq 1$.

\subhead I.1.2. The definition of quasideterminants\endsubhead  We start
from a definition for free skew-fields [GR1, GR2]. Let
$I,J$ be ordered sets consisting of $n$ elements.  Let $A=(a_{ij}), i\in
I, j\in J$ be a matrix with formal noncommuting entries $a_{ij}$.  Let us
define by induction $n^2$ rational expression $|A|_{pq}$ of variables
$a_{ij}, p\in I, q\in J$ over a free skew-field generated by $a_{ij}$'s.

We call these expressions the {\it
quasideterminants}.
\medskip
For $n=1$ we set $|A|_{ij} = a_{ij}$.
\medskip
For a matrix $A=(a_{ij}), i\in I, j\in J$ of order $n$ we denote
$A^{\alpha\beta}, \alpha\in I, \beta\in J$ the matrix of order $n-1$
constructed by deleting the row with the index $\alpha$ and the column
with the index $\beta$ in the matrix $A$.  Suppose that for the given
$p\in I, q\in J$ expressions $|A^{pq}|^{-1}_{ij}, i\in I, j\in J,
i\neq p, j\neq q$ are defined and set 
$$|A|_{pq} = a_{pq} -\sum
a_{pj}|A^{pq}|^{-1}_{ij} a_{iq}\tag 1.1.1
$$
Here the sum is taken over all $i\in I\smallsetminus\{p\},j\in
J\smallsetminus\{q\}$.
\definition{Definition 1.1.3}  We call the expression $|A|_{pq}$ the
quasideterminant of indices $p$ and $q$ of the matrix $A$.
\enddefinition 
The following definition of quasideterminants for matrices over general rings
was given in [GR].
\proclaim{Definition 1.1.4}  For a matrix A over a ring with unit 
the quasideterminant $|A|_{pq}$ is
defined if the matrix $A^{pq}$ is invertible.  In this case 
$$
|A|_{pq} = a_{pq} - \sum\Sb i\in I\smallsetminus\{p\}\\ j\in
J\smallsetminus\{q\}\endSb a_{pj}b_{ji} a_{iq},\tag 1.1.2
$$  
where $b_{ji}$ are the entries of the matrix $(A^{pq})^{-1}$.
\endproclaim
It is known [GR,GR1 ] that both definitions of a
quasideterminant coincide over free skew-fields.

\example{Examples}\endexample

1) For a matrix $a=(a_{ij}), i,j =1,2$
there exist four quasideterminants
$$\matrix
|A|_{11} = a_{11} - a_{12}\cdot a^{-1}_{22}\cdot a_{21},\quad
&|A|_{12}=a_{12}-a_{11}\cdot a^{-1}_{21}\cdot a_{22},\\
|A|_{21} = a_{21} - a_{22}\cdot a^{-1}_{12}\cdot a_{11},\quad &
|A|_{22}=a_{22}-a_{21}\cdot a^{-1}_{11}\cdot a_{12}.\endmatrix 
$$

2) For a matrix $A=(a_{ij}), i,j=1,2,3$ there exist 9 quasideterminants
but we will write here only 
$$\matrix
|A|_{11}=a_{11}-a_{12}(a_{22}-a_{23}a_{33}^{-1}a_{32})^{-1}a_{21} &-a_{12}(a_{32}-a_{33}\cdot
a_{23}^{-1} a_{22})^{-1} a_{31}\\
 \qquad -a_{13}(a_{23}-a_{22}a_{32}^{-1}a_{33})^{-1}a_{21} &-a_{13}(a_{33}-a_{32}\cdot
a^{-1}_{22}a_{23})^{-1}a_{31}\endmatrix
$$

3) If in formulas (1.1), or (1.2) variables $a_{ij}$ commute each
other then
$$
|A|_{pq} = (-1)^{p+q} \frac{\det A}{\det A^{pq}}.
$$


\remark{Remark} If any $a_{ij}$ is an invertible morphism $V_j\to V_i$ in
an additive category then the quasideterminant $|A|_{pq}$ is also a
morphism from the object $V_q$ to the object $V_p$.\endremark

One can use formulas (1.1) or (1.2) to define a quasideterminant
over any skew-field.  We recall here that according to Amistsur's
theory a rational function $f$ is defined over a ring if at least one
of the rational expressions for the function $f$ is defined over this
ring.  In particular, the quasideterminant $|A|_{pq}$ can be defined
if some of $(n-1)$-quasiminors of the matrix $A$ are not defined and
expression (1.1) makes no sense (but an equivalent expression does!).

The following theorem was conjectured by I. Gelfand and V. Retakh, and
proved by Reutenauer [R].
\proclaim{Theorem I.1.5}  Quasideterminants of the $n \times n$-matrix
with formal entries have the inversion height $n-1$ over the free skew-field
generated by its entries.
\endproclaim
{\centerline \bf {\S1.2 General properties of quasideterminants}} 
\subhead I.2.1. Two involutions\endsubhead
For a square matrix $A=(a_{ij})$ denote by $IA=A^{-1}$  inverse matrix,
and by $HA=(a^{-1}_{ji})$ its Hadamard inverse.  It is evident that if
$IA$, or $HA$ are defined then $I^2A=A$,respectively  $H^2A=A$. 

Let $A^{-1}=(b_{ij})$. According to [GR, GR1]
$b_{ij}=|A|^{-1}_{ji}$. This formula could be written in the following form. 
\proclaim{Theorem 1.2.1} For a square matrix A with formal entries
$$
HI(A)=(|A|_{ij})\tag 1.2.1
$$
\endproclaim
\subhead I.2.2. Homological relations [GR, GR1, GR2]\endsubhead  
The ratio of two
quasideterminants of the same square matrix is the ratio of two
rational functions of a smaller height.  It follows from the following
{\it homological} relations:
\proclaim{Theorem 1.2.3}

a) Row homological relations:
$$
-|A|_{ij} \cdot |A^{i\ell}|^{-1}_{sj} = |A|_{i\ell}\cdot
|A^{ij}|^{-1}_{s\ell}\qquad \forall s\neq i
$$

b) Column homological relations:
$$
-|A^{kj}|^{-1}_{it} \cdot |A|_{ij} = |A^{ij}|^{-1}_{kt}\cdot |A|_{kj}
\qquad \forall r\neq j
$$
\endproclaim
\subhead I.2.3 Heredity (special case)\endsubhead  Let $A=\pmatrix A_{11}&
A_{12}\\A_{21}& A_{22}\endpmatrix$ be a decomposition of $A=(a_{ij}),
i,j=1,\dots, n$ into a block-matrix.  Let $A_{11}$ be a $k\times
k$-matrix and the matrix $A_{22}$ be invertible.

\proclaim{Theorem 1.2.3} 
$|A|_{ij}= |A_{11}-A_{12} A_{22}^{-1}\cdot A_{21} |_{ij}
\text{ for } i,j=1.\dots, k.$  In other words, the quasideterminant
$|A|_{ij}$ of an $n\times n$ matrix can be computed in two steps: first,
considering the quasideterminant $|\tilde A|_{11} = A_{11} - A_{12}
\cdot A_{22}^{-1}\cdot A_{21}$of a $2\times 2$-matrix $\tilde
A=\pmatrix A_{11} &A_{12}\\
A_{21} & A_{22}\endpmatrix$; and second, considering the corresponding
quasideterminant of the $k\times k$-matrix $|\tilde A|_{11}$.
\endproclaim

\subhead I.2.4. General Heredity \endsubhead  Let $A=(a_{ij}),
i,j=1,\dots, n$ be a matrix and 
$$
A=\pmatrix A_{11} & \dots & A_{1s}\\
 {}&{} &{}\\
A_{s1} &\dots &A_{ss}\endpmatrix
$$
its block decomposition.  Denote by $\tilde A= (A_{ij}), i,j=1,\dots,
s$ the matrix with $A_{ij}$'s as entries.

Suppose that 
$$A_{pq} = \pmatrix
a_{k\ell} & \dots & a_{k,\ell+m}\\
{} &\dots &{}\\
a_{k+m,\ell} & \dots & a_{k+m,\ell+ m}\endpmatrix.
$$
is a square matrix and that $|\tilde A|_{pq}$ is defined.

\proclaim{Theorem 1.2.4}
$$
|A|_{k'\ell'}=\big\vert|\tilde A|_{pq}\big\vert_{k'\ell'} \quad\text{ for }
\matrix k'=k,\dots, k+m\\
\ell'=\ell,\dots, \ell +m\endmatrix.
$$
\endproclaim
\subhead I.2.5 Elementary properties of quasideterminants\endsubhead

\roster
\item"{i)}"  The quasideterminant $|A|_{pq}$ does not depend of the permutation
of rows and columns in the matrix $A$ if the $p$-$th$ row and the $q $-$th$
column are not changed;

\item"{ii)}"  {\it The multiplication of rows and columns.}  Let the
matrix $B$ be constructed from the matrix $A$ by multiplication of its
$i$-$th$ row by a scalar $\lambda$  from the left.  Then 
$$
 |B|_{kj}=\cases \lambda |A|_{ij} \qquad&\text{ if } k = i\\
|A|_{kj} \qquad&\text{ if } k \neq i \text{ and } \lambda \text{ is
invertible.}\endcases
$$

Let the matrix $C$ is constructed from the matrix A by multiplication
of its $j$-$th$ column by a scalar $\mu$ from the right.  Then
$$
|C|_{i\ell}=\cases |A|_{ij} \mu \qquad&\text{ if } \ell = j\\
|A|_{i\ell} \qquad&\text{ if } \ell \neq j \text{ and } \mu \text{ is
invertible.}\endcases
$$

\item"{iii)}" {\it The addition of rows and columns.} Let the matrix
$B$ is constructed by adding to some row of the matrix $A$ its
$k$-$th$ row multiplied by a scalar $\lambda$ from the left.  Then
$$
|A|_{ij} = |B|_{ij}, \quad  i=1, \dots k-1, k+1,\dots n,
 j=1, \dots, n.
$$
Let the matrix $C$ is constructed by addition to some column of
the matrix $A$ its $\ell$-$th$ column multiplied by a scalar $\lambda$ from
the right.  Then 
$$
|A|_{ij}= |C|_{ij} , \,  j=1,\dots ,\ell -1,\ell +1,\dots n, 
i= 1,\dots, n.
$$
More generally, let $X,Y$ be any $n\times n$-matrices such that the $i-th$
row of the matrix $X$ consists of the elements $\delta_{ik}$ and the
$j-th$ column of the matrix $Y$ consists of the elements $\delta_{j\ell}$(here
$\delta$ is the symbol of Kronecker).  Then 
$$
\matrix
&|XA|_{ij} = |A|_{ij} &i=k, j=1,\dots, n;\\
&|AY|_{ij} = |A|_{ij} &j=\ell, i= 1,\dots, n.
\endmatrix
$$ 
>From these identities we have 
\endroster

\proclaim{Proposition 1.2.5}  Let the quasideterminant $|A|_{ij}$ be
defined.  Then the following statements are equivalent.
\roster
\item"{i)}" $|A|_{ij} = 0$
\item"{ii)}"  the $i$-$th$ row of the matrix $A$ is a left linear
combination of the other rows of this matrix;
\item"{iii)}"  the $j$-$th$ column of the matrix $A$ is a right linear
combination of the other columns of this matrix.
\endroster
\endproclaim

There exists  a notion of linear dependence of rows and columns for
matrices over a skew-field if rows are multiplied by the elements of
this skew-field from the left and columns from the right.  So there
exist notions of the row rank and the column rank and these ranks are
equal [C, GR, GR1, GR2].  It also follows from Proposition 1.2.5 

By definition, a  $r$-quasiminor of a square matrix $A$ is a 
quasideterminant of its $r\times r$-submatrix.
\proclaim{Proposition 1.2.6}  The rank of the matrix $A$ over a skew
field is greater or equal to $r$ if and only if at least one of the
$r$-quasiminors of the matrix $A$ is defined and is not equal to zero.
\endproclaim
As it was pointed out in [KL], homological relations lead immediately
to the following analogue of the classical expansion of a
determinant by a row or a column.

\proclaim{Proposition 1.2.7} For any $k\neq p$ and any $\ell \neq q$
$$ |A|_{pq}=a_{pq}-\sum\Sb j\neq q\endSb a_{pj}(|A^{pq}|_{kj})^{-1}
|A^{pj}|_{kq}$$

$$|A|_{pq}=a_{pq}-\sum\Sb i\neq p\endSb |A^{iq}|_{pi}(|A^{pq}|_{i\ell})^{-1}
a_{iq}$$
if all terms in these expressions are defined.
\endproclaim

{\centerline {\bf\S 1.3 Basic identities}}

\subhead I.3.1 Basic identities\endsubhead
The following noncommutative analogue of Sylvester's identity is
closely related with the heredity property.  Let $A=(A_{ij}),
i,j=1,\dots, n$ be a matrix over a noncommutative skew-field.  Suppose
that its submatrix $A_0=(a_{ij}), i,j=1, \dots. k$ is invertible.  For
$p,q=k+1,\dots, n$ set
$$
b_{pq} = \vmatrix &A_0 & {} & a_{1q}\\
&{} & {} & \vdots\\
&{} &{} &a_{kq}\\
&{} &{} &\dots\\
&a_{p1}\dots a_{pk}&\dots &a_{pq}\endvmatrix_{pq} ;
$$
$$
B=(b_{pq}), p,q = k+1,\dots, n.
$$
We call the submatrix $A_0$ {\it a pivot} for matrix $B$.
\proclaim{Theorem 1.3.1} For $i,j = k+1,\dots, n$
$$
|A|_{ij} = |B|_{ij}
$$
\endproclaim
For example, a quasideterminant of an $n\times n$-matrix $A$ is
equal to the corrseponding quasideterminant of the $2\times 2$-matrix 
consisting of $(n-1)$- quasiminors,
or to the quasideterminant of $(n-1)\times(n-1)$-matrix consisting of
$2$-quasiminors of the matrix $A$. One can use any of these
procedures for a definition of a quasideterminants.

In a special case when $A_0=(a_{ij}), i,j=2,\dots,
n-1$, Theorem 1.3.1 is an analogue of a well-known commutative indentity
which goes under name of ``Lewis Carrol identity''. 

\proclaim{Corollary 1.3.2} (Sylvester):  Let $A=(a_{ij}), i,j=1,\dots,
n$ be a matrix over a commutative ring.  Suppose that its submatrix
$A_0=(a_{ij}), i,j=1,\dots, k$ is invertible.  For $p,q = k+1,\dots,
n$ set
$$
\aligned
\tilde b_{pq} &= \det\pmatrix &A_0 & a_{1q}\\
&{}&\vdots\\
&{}&a_{kq}\\
&{}&\dots \\
&a_{p1}\dots a_{pk} &\dots a_{pq}\endpmatrix ,\\
\tilde B&=(\tilde b_{pq}),\ p,q= k+1,\dots, n.
\endaligned
$$
  
Then 
$$
\det A={\det\tilde B\over (\det A_0)^{n-k-1}}.
$$
\endproclaim
\subhead I.3.2. A noncommutative version of Jacobi theorem (involution
for quasiminors)\endsubhead  For a matrix $A=(a_{ij}), i\in I, j\in
J, P\subset I, Q\subset J$ denote by $A_{PQ}$ the submatrix:
$$
A_{PQ} = (a_{\alpha\beta}),\alpha\in P,\beta\in Q.
$$

Let $|I|= |J|$ and $B=A^{-1} = (b_{rs})$.  Suppose that $|P| = |Q|$.

\proclaim{Theorem 1.3.3}  Let $k\notin P, \ell\notin Q$.  Then
$$
|A_{P\cup\{k\}, Q\cup\{\ell\}}|_{k\ell}\cdot |B_{I\setminus Q,
J\setminus P}|_{\ell k} = 1.$$
\endproclaim
Set $P=I\smallsetminus\{k\}, Q= J\smallsetminus\{\ell\}$.  Then
this theorem leads to the already mentioned identity
$$
|A|_{k\ell} \cdot b_{\ell k} = 1.
$$

\subhead I.3.3. A generalization of the homological relations\endsubhead

For a matrix $A=(a_{ij}), i\in I, j\in J$ and two subsets $L\subset I,
M\subset J$ denote by $A^{L,M}$ the submatrix of the matrix $A$
obtained from $A$ by deleting the rows with the indexes $\ell \in L$
and the columns with the indexes $m\in M$.

Let $A$ be a square matrix, $L= (\ell_1,\dots, \ell_k), M=(m_0,\dots,m_k).$
 Set 
$M_i = M\smallsetminus\{m_i\}, i=0,\dots, k$.
\proclaim{Theorem 1.3.4} [GR1, GR2]  For $p\notin L$
$$
\sum^k_{i=0}
|A^{L,M_i}|_{pm_i}\cdot|A|^{-1}_{\ell m_i}= \delta_{p\ell},
$$
$$
\sum^k_{i=0}
|A|^{-1}_{m_i\ell }\cdot|A^{M_i,L}|_{m_ip} = \delta_{\ell p}.
$$
\endproclaim

\subhead I.3.4 Multiplicative properties of quasideterminants\endsubhead

Let $X =(x_{pq}), Y=(y_{rs})$ be $n\times n $-matrices.

\proclaim{Theorem 1.3.5}
$$|XY|^{-1}_{ij} = \sum^n_{p=1}|Y|^{-1}_{pj} |X|^{-1}_{pi}.
$$
\endproclaim
\medskip
{\centerline {\bf \S 1.4 Noncommutative linear algebra}}
\medskip
Here we recall some results from [GR], [GR1], [GR2].

\subhead I.4.1. Solutions of systems of linear equations\endsubhead

\proclaim{Theorem 1.4.1} If the system
$$
\cases
a_{11} x_1+\dots &+ a_{1n} x_n = \xi_1\\
\dots \\
a_{n1} x_1+\dots &+ a_{nn} x_n = \xi_n
\endcases
$$
is given and for the matrix $A=(a_{ij})$ the quasideterminants of $A$
are defined and invertible, then 
$$
x_i = \sum^n_{j=1} |A|^{-1}_{ji} \xi_j.
$$
for $i=1,\dots, n$.
\endproclaim

\subhead I.4.2 Cramer rule\endsubhead  Let $A_{\ell}(\xi)$ be the $
n\times n$-matrix constructed by replacing in the matrix $A$ its
$\ell-th$ column by the column $(\xi_1,\dots, \xi_n)$.

\proclaim{Theorem 1.4.2}  If the quasideterminants $|A|_{ij}$ and
$|A_j(\xi)|_{ij}$ are defined then
$$
|A|_{ij} x_j = |A_j(\xi) |_{ij}.
$$
\endproclaim
\subhead I.4.3. Cayley - Hamilton theorem\endsubhead  Let $t$ be a
formal variable and $A=(a_{ij}), i,j=1,\dots, n$ be $a$ matrix over a
skew field.  We do not suppose that $t$ commutes with variables
$a_{ij}, i,j=1,\dots, n$.  Let $\Bbb I_n$ be $n\times n$-unit matrix
and $f_{ij}(t) = |t\cdot \Bbb I_n- A|_{ij}$.  Consider the expression
$\tilde f_{ij}(t)$ as a matrix function of a matrix variable
changing in  $f_{ij}(t)$
elements $a_{ij}$ by matrices $ \tilde a_{ij} = a_{ij}\cdot \Bbb I_n$.
For example, if $n=2$ then $\tilde f_{11} (t) = t-\tilde a_{11} -
\tilde a_{12} (t-\tilde a_{22})^{-1} \tilde a_{21}$. 

The following
theorem was formulated in [GR1], [GR2]. An independent proof
(similar to the original one) was given by C. Reutenauer in his
Montreal lectures.

\proclaim{Theorem 1.4.3}  For all $i,j=1,\dots, n$
$$
f_{ij}(A) = 0.
$$
\endproclaim

\head II. Noncommutative Pl\"ucker and Flag Coordiniates\endhead
\medskip
{\centerline {\bf \S2.1 Definition and basic properties of 
quasi-Pl\"ucker coordinates}}

\subhead II.1.1.  Definition of left quasi-Pl\"ucker coordinates\endsubhead
  Let $A=(a_{pq})$,$\ p=1, \dots, k,\ q=1, \dots, n$,$ k< n$ be a matrix over
a skew-field $F$.  
Fix 

$1\leq i, j, i_1,\dots, i_{k-1}\leq n$ such
that $i\notin I=\{ i_1,\dots, i_{k-1}\}$.  For $1\leq s\leq k$ set
$$
p^{I}_{ij}(A) =\vmatrix &a_{1i}a_{1i_1}&\dots & a_{1i_{k-1}}\\
&{}  &\dots & {}\\
&a_{ki}a_{ki_1} & \dots &a_{ki_{k-1}}\endvmatrix^{-1}_{si} \cdot
\vmatrix &a_{1j} a_{1i_1} &\dots &a_{1i_{k-1}}\\
& {}  &\dots  &{}\\
&a_{kj} a_{ki_1} &\dots &a_{ki_{k-1}}\endvmatrix_{sj}.
$$

\proclaim{Proposition 2.1.1}  

$i) p_{ij}^I (A)$ does not depend of $s$;

$ii)  p^I_{ij}(gA)=p^I_{ij}(A)$ for any invertible $k\times k$ matrix $g$.
\endproclaim

\demo{Proof} Use homological relations.\enddemo

\definition{Definition 2.1.2} We call $p^I_{ij}(A)$ left
quasi-Pl\"ucker coordinates of the matrix $A$.  
\enddefinition
One may consider $p^I_{ij}$ as quasi-Pl\"ucker coordinates of a
$k$-dimensional left subspace generated in the left $F$-linear space $F^n$ by
rows of the matrix $A$.

In commutative case $p^I_{ij} = {p_{jI}\over p_{iI}}$, where
$p_{\alpha_1\dots \alpha_k}$ is the standard Pl\"ucker
coordinate.

\subhead II.1.2 Generating identities \endsubhead  The following
properties of $p^I_{ij}$'s trivially follow from their definition:
\roster 
\item"{i)}" $p^I_{ij}$ does not depend of reordering of elements in $I$;
\item"{ii)}" $p^I_{ij} = 0 \text{ if } j\in I;$
\item"{iii)}"  $p^I_{ii} = 1 \text{ and } p^I_{ij} \cdot p^I_{jk} = p^I_{ik}$.
\endroster

\proclaim{Proposition 2.1.3} \text {(Skew-Symmetry)}.  Let $N$ be a set
of indices, $|N|=k+1, i,j,m\in N$.  Then 
$$
p^{N\smallsetminus\{i,j\}}_{ij} \cdot p^{N\smallsetminus\{j,m\}}_{jm}\cdot
p^{N\setminus\{m,i\}}_{mi} =- 1.
$$
  
\endproclaim

\proclaim{Proposition 2.1.4} \text {(Pl\"ucker relations)}.  
Fix 
$M=(m_1,\dots, m_{k-1})$, $L=(\ell_1,\dots, \ell_k)$. Let $i\notin M$. Then
$$
\sum_{j\in L} p^M_{ij} \cdot p_{ji}^{L\setminus\{j\}}=1.
$$
\endproclaim

\example{Examples} Suppose that $k=2$.

1) From Proposition 2.1.3 it follows that 
$$
p^{\{\ell\}}_{ij}\cdot p^{\{i\}}_{j\ell}\cdot p^{\{j\}}_{\ell i}=-1.
$$

In the commutative case one has $p^{\{\ell\}}_{ij} = {p_{j\ell}\over
p_{i\ell}}$ and $p_{ij} =- p_{ji}$, so this identity is valid.

2)  From Proposition 2.1.4 it follows that for any $i,j,\ell, m$
$$
p^{\{\ell\}}_{ij} \cdot p_{ji}^{\{m\}} + p_{im}^{\{\ell\}}\cdot
p^{\{j\}}_{mi} =1.
$$

In commutative case this identity implies the standard identity
$$
p_{ij}\cdot p_{\ell m} - p_{i\ell}\cdot p_{jm} + p_{im}\cdot p_{\ell
j} =0.
$$
\endexample
\remark{Remark} In fact the products $p^{\{\ell}\}_{ij}p^{\{m}\}_{ji}$ 
(which are equal to 
${p_{j\ell}\over p_{i\ell}}\cdot {p_{im}\over p_{jm}}$ 
in the commutative case) are 
noncommutative cross-ratio's.
\endremark

We will consider different  specializations of the relations of
quasi-Pl\"ucker coordinates in a separate paper.

To prove propositions 2.1.3 and 2.1.4 we need the following lemma.
Let $A=(a_{ij}), i=1,\dots ,k, j=1,\dots , n, k< n$ be a matrix over 
a skew-field. Denote by $A_{j_1,\dots ,j_k}$ its submatrix $(a_{ij}),
i=1,\dots ,k, j=j_1,\dots ,j_k$. Consider $n\times n$-matrix

$$X=\pmatrix A_{1\dots k}&A_{k+1\dots n}\\0&E_{n-k}\endpmatrix ,$$
where $E_m$ is a $m\times m$-unit matrix.

\proclaim{Lemma 2.1.5} Let $j< k< i$. If 
$p^{1\dots \hat {j}\dots k}_{ij}(A)$ is defined then $|X|_{ij}$ is
defined and

$$|X|_{ij}=-p^{1\dots \hat {j}\dots k}(A) \tag 2.1 $$
\endproclaim

{\bf Proof.} We have to prove

$$|X|_{ij}=-|A_{1\dots \hat {j}\dots k i}|^{-1}_{si}
\cdot |A_{1\dots k}|_{sj} \tag 2.2 $$
if the right-hand side is defined. We will prove it by induction over
$\ell =n-k$. Suppose that formula (2.2) is valid for $l=m$ and prove it
for $\ell =m+1$. Without loss of generality we may set $j=1, i=k+1$.
By Proposition 1.2.7
$$|X|_{k+1,1}=-|X^{k+1,1}|^{-1}_{s,k+1}\cdot |X^{k+1,k+1}|_{s1}$$

for an appropriate $1\leq s\leq k$. Here
$$X^{k+1,1}=
\pmatrix A_{2\dots k+1}&A_{k+2\dots n}\\0&E_{n-k-1}\endpmatrix ,$$

$$X^{k+1,k+1}=
\pmatrix A_{1\dots k}&A_{k+2\dots n}\\0&E_{n-k-1}\endpmatrix .$$

By the induction assumption
$$|X^{k+1,1}|_{s,k+1}=
-|A_{23\dots k k+2}|^{-1}_{s,k+2}\cdot |A_{23\dots k+1}|_{s,k+1}, $$

$$|X^{k+1,k+1}|_{s1}=
-|A_{23\dots k k+2}|^{-1}_{s,k+2}\cdot |A_{1\dots k}|_{s1} $$

and $|X|_{k+1,1}=-p^{23\dots k}_{k+1,1}$. $\square $

To prove Proposition 2.1.4 we apply Theorem 1.3.4 ii) to the matrix

$$X=\pmatrix A_{1\dots k}&A_{k+1\dots n}\\0&E_{n-k}\endpmatrix $$

for $M=(k+1,\dots ,n)$ and any $L$ such that $|L|=n-k-1$. By Lemma 2.1.5

$|X|_{m_i\ell }=-p^{1\dots \hat \ell \dots k}(A)$,
$|X^{M_i,L}|_{m_i q}=-p^{{1\dots n}\setminus L}_{m_i q}(A)$, and
Proposition 2.1.4 follows from Theorem 3.4.ii) $\square $.

To prove Proposition 2.1.3 one have consider matrix $X$ for $n=k+1$ and
use homological relations.

\proclaim{Theorem 2.1.6}  Let $A =(a_{ij}), i=1, \dots, k, j=1,\dots, n$
be a matrix with formal entries, and $f (a_{ij})$ be a rational function
over a free skew-field $F$ generated by $(a_{ij}$.  Let $f$ be invariant 
for all transformations:
$$
A\to gA
$$
where $g$ is an invertible $k\times k$-matrix over $f$.  Then $f$ is a rational
function of quasi-Pl\"ucker coordinates.
\endproclaim
{\bf Proof}. Let $b_{ij}=a_{ij}$ for $i,j=1,\dots k$. Consider matrix
$B=(b_{ij})$. Then $B^{-1}=(|B|^{-1}_{ji})$. 

Set $C=(c_{ij})=B^{-1}A$. Then 
$$
c_{ij}=\cases \delta _{ij}\qquad& j\leq k \\
p^{1\dots \hat {i}\dots k}_{ij}(A)\qquad& j>k
\endcases $$

By invariance we have $f(A)=f(1,C)$. $\square $
              
\subhead II.1.4. Right quasi-Pl\"ucker coordinates\endsubhead
 Consider a matrix $B=(b_{ij}), i=1,\dots, n;
j=1,\dots, k, k< n$ over a skew-field $F$.  Fix $1\leq i, j,
i_1,\dots, i_{k-1} \leq n$ such that $j\notin I = (i_1,\dots,
i_{k-1})$.  Fix also $1\leq t\leq k$ and set
$$
r^I_{ij}(B) = \vmatrix &b_{i1} &\dots & b_{ik}\\
&b_{i_11} & \dots & b_{i_1k}\\
&{} &\dots &{}\\
&b_{i_{k-1}1} &\dots &b_{i_{k-1}k}\endvmatrix_{it} \cdot 
\vmatrix &b_{j1} & \dots &b_{jk}\\
&b_{i_11} &\dots &b_{i_1 k}\\
&{} &\dots &{} \\
&b_{i_{k-1} 1} &\dots &b_{i_{k-1}k}\endvmatrix^{-1}_{jt} .
$$

\proclaim{Proposition 2.1.7} 

i) $r^I_{ij}(B)$ does not depend of
$t$; 

ii) $r^I_{ij}(Bg) = r^I_{ij}(B)$ for any invertible
$k\times k$-matrix 
$ g\text{ over } F$
\endproclaim
\proclaim{Definition 2.1.8} We call $r^I_{ij}(B)$ right quasi-Pl\"ucker
coordinates of the matrix $B$.
\endproclaim
Proposition 2.1.7 ii) shows that $r^I_{ij}(B)$ are quasi-Pl\"ucker
coordinates of the right $k$-dimensional subspace generated by columns
of $B$ in the right vector space $F^n$.

\subhead II.1.5  Generating identities for right quasi-Pl\"ucker
coordinates.\endsubhead  Generating identities for $r^I_{ij}$'s  
are dual to generating 
identities for the left quasi-Pl\"ucker coordinates
$p^I_{ij}$.

\subhead II.1.6 The duality between quasi-Pl\"ucker
coordinates\endsubhead

Let $A=(a_{ij}), i=1,\dots, k; j=1,\dots, n;$ and $B=(b_{rs}), r=1,\dots,
n; s=1,\dots, n-k$.  Suppose that $AB=0$.  It is equivalent to the
statement that the subspace generated by rows of $A$ in the left
linear space is orthogonal to the subspace generated by columns of $B$
in the right linear space.  Fix $1\leq i, j\leq n, I\subset [1,n],
i\notin I, |I| = k-1$.  Set $J=[1,n]\setminus I\setminus \{i,j\}$

\proclaim{Theorem 2.1.9}
$$
p^I_{ij}(A) + r^J_{ij} (B) = 0.
$$\endproclaim
 
\subhead II.1.7 Quasi-Pl\"ucker coordinates for $k\times n$-matrices for
different k\endsubhead
 
Let  $A=(a_{\alpha\beta}), \alpha=1,\dots, k; \beta = 1, \dots, n$ be any
$k\times n$ - matrix over a noncommutative skew-field and $A'$ be a
$(k-1)\times n$-submatrix of $A$.  Fix  $1\leq i, j, m, j_1,\dots,
j_{k-2} \leq n$  such that $i,m\notin J=\{j_1, \dots, j_{k-2}\}$ and
$i\neq m$.

\proclaim{Proposition 2.1.10} 
$$
p^J_{ij}(A')= p_{ij}^{J\cup\{m\}}(A) + p^J_{im}(A')\cdot
p_{mj}^{J\cup\{i\}}(A).
$$
\endproclaim

{\centerline {\bf\S 2.2 Applications of quasi-Pl\"ucker coordinates}} 

\subhead II.2.1  Row and column expansion of a
quasideterminant\endsubhead

Some of the results obtained in [GR], [GR1], [GR2] and partially
listed in Chapter I can be rewritten in terms of quasi-Pl\"ucker
coordinates. 
Let $A=(a_{ij}), i, j=1, \dots, n$  be a matrix over a skew-field.

Fix $1\leq \alpha,\beta\leq n $. Let $B=(a_{ij}, i\neq \alpha ),
C=(a_{ij}, j\neq \beta )$ be its $(n-1)\times n$ and $n\times (n-1)$
submatrices respectively.  
For any $j\neq \beta$ and $i\neq
\alpha$ set:
$$
\aligned
p_{j\beta} &= p^{1\dots\hat j \dots \hat \beta\dots n}_{j\beta}
(B),\\
r_{\alpha i}&= r^{1\dots\hat\alpha\dots\hat i \dots n}_{\alpha i}
(C).\endaligned
$$

\proclaim{Proposition 2.2.1}

i)  $|A|_{\alpha \beta} = a_{\alpha\beta} - \sum_{j\neq \beta}
a_{\alpha j} p_{j\beta}$,

ii)  $|A|_{\alpha\beta} = a_{\alpha\beta} - \sum_{i\neq \alpha}
r_{\alpha i} a_{i\beta}$

if terms in the right-hand sides are defined.
\endproclaim

\subhead II.2.2. Homological relations\endsubhead

\proclaim{Proposition 2.2.2}  In notations of II.2.1

i) $|A|^{-1}_{ij} \cdot |A|_{i\ell} =- p_{j\ell}$  \qquad (row
relations)

ii)  $|A|_{ij} \cdot |A|^{-1}_{kj} =- r_{ik}$      \qquad (column
relations).

\endproclaim
\proclaim{Corollary 2.2.3}  In notations of II.2.1, let $(i_1,\dots, i_s),
(j_1,\dots, j_t)$ be sequences of indices such that $i\neq i_1, i_1\neq
i_2,\dots, i_{s-1}\neq i_s; j \neq j_1, j_1\neq j_2,\dots, j_{t-1}\neq
j_t$.  

Then 
$$ 
|A|_{i_sj_t} = p_{i_si_{s-1}}\dots p_{i_2i_1} p_{i_1i}\cdot |A|_{ij}
\cdot r_{jj_1}r_{j_1j_2}\dots r_{j_{t-1}j_t}.
$$
\endproclaim
\example{Example} \endexample  For a matrix $A=\pmatrix a_{11}&
a_{12}\\ a_{21} & a_{22}\endpmatrix$
$$
\aligned
&a_{21} \cdot a^{-1}_{11}\cdot |A|_{11} \cdot a^{-1}_{22}\cdot a_{22}
= |A|_{22},\\
&a_{12} a_{22}^{-1}a_{21} a^{-1}_{11} |A|_{11} a^{-1}_{21} a_{22}
a^{-1}_{12} a_{11} = |A|_{11}\endaligned
$$
\subhead II.2.3 Matrix multiplication\endsubhead The following formula
has been already used in the proof of Theorem 2.1.6.
Let $A= (a_{ij}), i=1, \dots, n; j=1,\dots, m, n< m; B=(a_{ij}), i =
1,\dots, n, j=1, \dots, n; C=(a_{ik}), i=1,\dots, n; k=n+1, \dots,
m.$

\proclaim{Proposition 2.2.4} Let matrix $B$ is invertible. Then
$p^{1\dots \hat i\dots n}_{ik}(A)$ are defined for $i=1,\dots n,
k=n+1,\dots m$ and 
$$
B^{-1}C=(p^{1\dots\hat i\dots n}_{ik}(A)), i=1,\dots, n; k= n+1,
\dots, m.
$$
\endproclaim

\subhead II.2.5 Quasideterminant of a product\endsubhead
Let $A=(a_{ij}), B=(b_{ij}), 
i,j=1,\dots n$ be matrices over a skew-field. Fix $1\leq k\leq n$.
Consider $(n-1)\times n$-matrix $A'=(a_{ij})$, $i\neq k $, 
and $n\times (n-1)$-matrix $B''=(b_{ij})$,  $j\neq k$.

\proclaim{Proposition 2.2.5}
$$|B|_{kk}\cdot |AB|^{-1}_{kk}|A|_{kk} = 
1+\sum_{\alpha \neq k} r_{k\alpha}\cdot
p_{\alpha k},$$

where
$$
\aligned
r_{k\alpha} &= r^{1\dots\hat \alpha\dots n}_{k\alpha}(B'')-
\text{ right quasi-Pl\"ucker coordinates}\\
p_{\alpha k} &= p_{\alpha k}^{1\dots\hat\alpha\dots n}(A')-
\text{ left quasi-Pl\"ucker coordinate}
\endaligned
$$
if all the expressions are defined.
\endproclaim
The proof follows from the multiplicative property of quasideterminants
and Proposition 2.2.2.
\subhead II.2.6. Gauss decomposition \endsubhead Consider matrix
$A=(a_{ij}),i,j=1,\dots,n$ over a skew-field. Let $A_k=(a_{ij}),
i,j =k,\dots n$, $B_k=(a_{ij}), i=1,\dots n, j=k, \dots n$, and
$C=(a_{ij}), i=k,\dots n, j=1,\dots n$ be its submatrices of sizes
$(n-k+1)\times (n-k+1)$, $n\times (n-k+1)$, and $(n-k+1)\times n$
respectively.  

Suppose that quasideterminants
$$
y_k=|A_k|_{kk},\  k=1,\dots, n
$$
are defined and invertible. From [GR1], [GR2] it follows
\proclaim{Theorem 2.2.5}
$$
A=\pmatrix 1 &{}&x_{\alpha \beta}\\
{}&\ddots&{}\\
0&{}&1\endpmatrix
\pmatrix y_1&{}&0\\
{}&\ddots&{}\\
0&{}&y_n\endpmatrix\pmatrix1&{}&0\\
{}&\ddots&{}\\
z_{\beta \alpha}&{}&1\endpmatrix,
$$
where
$$
\aligned
x_{\alpha \beta} &=r^{\beta+1\dots n}_{\alpha\beta}(B_{\beta}),
\ 1\leq \alpha <\beta\leq n\\
z_{\beta\alpha} &=p^{\beta+1\dots n}_{\beta\alpha} ( C_{\beta}).
\ 1\leq \alpha <\beta\leq n\\
\endaligned
$$
\endproclaim

\subhead II.2.7 Flag coordinates \endsubhead  We remind now a
definition of noncommutative flag coordinates introduced in [GR1],
[GR2].
Let $A=(a_{ij}), i=1,\dots, k;
j=1,\dots,n$ be a matrix over a skew-field $R$.  Let
$\Cal F=(F_1\subset F_2\subset\dots\subset F_k)$ be a flag in the left
vector space $R^n, F_p$ is generated by the first $p$ rows of
$A$. Put
$$
f_{j_1 \dots j_k} (\Cal F) =\vmatrix a_{1j_1}& \dots  &a_{1j_k}\\
{}&\dots &{}\\
a_{kj_1}&\dots &a_{kj_k}\endvmatrix_{kj_1}.
$$

\proclaim{Proposition} [GR1], [GR2]
Functions $f_{j_1\dots j_m}(\Cal F)$ do not change under
multiplication of $A$ by upper triangular matrix with units as
its diagonal entries from the left.
\endproclaim
 In [GR1], [GR2] functions $f_{j_1\dots j_m}(\Cal F)$ 
were called flag coordinates of $\Cal F$. 

It is easy to see that
$$
p_{ij}^{i_1\dots i_{k-1}}
(A)=\left(f_{ii_1\dots i_{k-1}}(\Cal F)\right)^{-1}\cdot f_{ji_1\dots
i_{k-1}}(\Cal F).
$$
In [GR1],[GR2] 
the relations between flag
coordinates for $k\times n$- and $(k-1)\times n)$-matrices   
were considered. Our Proposition 2.1.9
may be deduced from these relations.

\head III. Noncommutative Symmetric Functions, Bezout Theorem, and Vieta
Theorem
\endhead
\S 3.1 Bezout and Vieta Theorems
\subhead III.1.1 Vandermonde quasideterminants\endsubhead
In this section we will discuss analogues of Bezout and 
Vieta Theorem for polynomials over a noncommutative skew-field. We will
prove Vieta Theorem which was formulated in [GR3] using our noncommutative
form of Sylvester identity. The first prove of this theorem using
differential operators has appeared in [EGR]. Another noncommutative
version of Vieta theorem based on notions of traces and determinants
was given by A. Connes and A. Schwarz in [CS].
 
\medskip
Let $x_1, x_2, \dots, x_k$ be a set of elements of a skew-field $F$.
For $k>1$ define a {\it Vandermonde} quasideterminant
$$
V(x_1,\dots,x_k)=\vmatrix x_1^{k-1} & \dots &x_k^{k-1}\\
{} &\dots &{}\\
x_1&\dots &x_k\\
1&\dots & 1\endvmatrix_{1k}.
$$

We say that a sequence $x_1,\dots, x_n\in F$ is {\it independent} of
for any $k=2, \dots, n$ quasideterminants $V(x_1,\dots, x_k)$ are
defined and invertible.  For independent sequences $x_1, \dots, x_n$
and $x_1,\dots , x_{n-1}$, where $x_1,\dots ,x_n, z\in F$
set 
$$
\aligned
y_1 &=x_1, \, z_1 = z\\
y_k &=V(x_1,\dots,x_k) x_k V^{-1}(x_1,\dots, x_k), k\geq 2\\
z_k &=V(x_1,\dots,x_{k-1}, z) z V^{-1}(x_1,\dots, x_{k-1},z), k\geq 2.
\endaligned
$$

In the commutative case $y_k = x_k$, and $z_k = z$ for $k=1,\dots, n$.

\subhead III.1.2 Bezout and Vieta decomposition\endsubhead   
We will prove here
\proclaim{Theorem 3.1.1} (Bezout decomposition of a Vandermonde
quasideterminant).  Suppose that sequences $x_1,\dots, x_n$ and
$x_1,\dots , x_{n-1}, z$ are independent. then
$$
V(x_1,\dots, x_n, z) = (z_n-y_n)(z_{n-1}-y_{n-1})\dots(z_1-y_1).\tag 3.1
$$
\endproclaim
Note that 
if $z$ commutes with $x_i, i=1,\dots n$ then
$$
V(x_1,\dots , x_n, z)= (z-y_n)(z-y_{n-1})\dots (z-y_1).
$$
\proclaim{Theorem 3.1.2} (Vieta decomposition of a Vandermonde
quasideterminant).  For a independent sequence $x_1,\dots, x_n, z$
$$
V(x_1,\dots, x_n, z) = z^n + a, z^{n-1} + \dots + a_{n-1} z+ a_n,\tag 3.2
$$
where 
$$
a_k=(-1)^k\sum_{1\leq i_1<i_2<\dots<i_k\leq n} y_{i_k} 
\cdot y_{i_{k-1}}\dots y_1 .\tag 3.3
$$

In particular
$$
\aligned
a_1&= - (y_1+\dots+ y_n),\\
a_2&= \sum_{1\leq i<j\leq n} y_jy_i,\\
 & \dots\\
a_n&=(-1)^ny_n\dots y_1.
\endaligned
$$
\endproclaim
Theorem 3.1.2 follows from Theorem 3.1.1 by induction.

\subhead III.1.3  Proof of Theorem 3.1.2\endsubhead  
For $n=1$ one has $V(x_1,z)= z-x_1$ and
formulas (3.2 ) are valid.  Suppose that these formulas are valid for
$m=n-1$.  We have to prove that formula (3.2) can be written in the
form 
$$
V(x_1,\dots, x_n, z) = (z_n-y_n) V(x_1,\dots, x_{n-1}, z) = 
$$
$$
=V(x_1,\dots, x_{n-1}, z)\cdot z- y_n\cdot V (x_1, \dots, x_{n-1}, z).
$$

By induction,
$$
V(x_1,\dots, x_{n-1},z)= z^{n-1} + b_1 z^{n-2} + \dots + b_{n-1},
$$
where
$$
\aligned
b_1 &=-(y_1+\dots+y_{n-1}),\\
{} &\dots\\
 b_{n-1}&=(-1)^n y_{n-1}\cdot\dots\cdot y_1.
\endaligned
$$

So, 
$$
V(x_1,\dots,x_n, z) = z^n+(b_1-y_n) z^{n-1}+(b_2-y_nb_1) z^{n-2} +
\dots -y_nb_n=$$
$$= z^n+a_1 z^{n-1} + \dots + a_n,
$$
where $a_1,\dots, a_n$ are given by formulas (3.3).\qed

\subhead III.1.4 \endsubhead  To prove Theorem
3.1.1 
we need the following 
\proclaim{Lemma 3.1.3}
$$
V(x_1,\dots, x_n, z) = V(\hat x_2,\dots, \hat x_n, \hat z)(z-x_1),
$$
where 
$$
\aligned
\hat x_k &=(x_k - x_1) x_k(x_k-x_1)^{-1}, k=2,\dots, n\\
\hat z &=(z-x_1) z(z-x_1)^{-1}.
\endaligned
$$
\endproclaim  
\demo{Proof} By definition
$$
V(x_1,\dots, x_n, z) =\vmatrix &x_1^n &x^n_2 &\dots &z^n\\
&x_1^{n-1} &x_2^{n-1} &\dots &z^{n-1}\\
 & {}&{} &\dots  &{}\\
&x_1 &x_2 &\dots &z\\
&1 &1 &\dots &1\endvmatrix_{1 n+1}.
$$
Multiply $k$-$th$ row by $x_1$ from the left and subtract it from
$(k-1)$-$th$ row for $k=2,\dots, n$.  The quasideterminant will not change,
and so 
$$
V(x_1,\dots, x_n z) = \vmatrix
&0 &x^n_2-x_1x_2^{n-1} &\dots &z^n-x_z^{n-1}\\
&0 &x_2^{n-1} - x_1x_2^{n-2} &{} &z^{n-1}-x_1 z^{n-2}\\
&\vdots &\vdots  &{} &\vdots\\
&0 &x_2-x_1 &{} &z-x_1\\
&1 & 1 &{} & 1\endvmatrix_{1 n+1} =
$$
$$
=\vmatrix &0 &(x_2-x_1)x_2^{n-1} &\dots &(z-x_1)z^{n-1}\\
&\vdots&\vdots &{} &\vdots\\
&0 & x_2-x_1& {} &z- x_1\\
&1  &1 &{} & 1\endvmatrix_{1 n +1}.
$$
Applying Sylvester theorem to the last quasideterminant with the
element of index $(n1)$ as a pivot one has 
$$
V(x_1,\dots, x_n, z) = \vmatrix 
&(x_2-x_1)x_2^{n-1}&\dots &(z-x_1)z^{n-1}\\
&\dots &{} &\dots\\
&x_2-x_1 &{} &z-x_1\endvmatrix_{1n} .
$$
Multiply $k$-$th$ column by $(x_{k+1} - x_1)$ for $k=1,\dots, n-1$ and
the last column by $(z-x_1)$.  According to elementary properties of
quasideterminants,
$$
V(x_1,\dots, x_n, z)=
$$
$$
=\vmatrix 
&(x_2-x_1)x_2^{n-1}(x_2-x_1)^{-1} &\dots &(z-x_1)z^{n-1}(z-x_1)^{-1}\\
&\vdots  &{} &\vdots\\
&1 & {} & 1\endvmatrix _{1n} (z-x_1) = 
$$
$$
=\vmatrix&\hat x_2  &\dots &\hat z^{n-1}\\
&\vdots  & {} &\vdots\\
&1   &\dots &1\endvmatrix _{1n}\cdot (z-x_1) = V(\hat x_2,\dots, \hat x_n,
\hat z)\cdot (z-x_1).$$
\hfill\qed
\enddemo

\subhead III.1.5 Proof of Theorem 3.1.1\endsubhead
Use induction.  By Lemma 3.1.3 
Theorem 3.1.1 is valid for $n=2$.  Also by Lemma 3.3.1
$$
V(x_1,\dots, x_n, z) = V(\hat x_2,\dots, \hat x_n, \hat z)(z-x_1).\tag 3.4
$$
Suppose that our theorem is valid for $m=n-1$.

Then
$$
V(\hat x_2,\dots,\hat x_n,\hat z) = (z'_n - y'_n)\dots (z'_2- y'_2),
$$
where 
$$
\aligned
z'_2 &=\hat z, \\ 
y'_2 &= \hat x_2, \\
z'_k &=V(\hat x_2,\dots, \hat x_{k-1},\hat z) \hat z V^{-1}(\hat
x_2,\dots,\hat x_{k-1},\hat z),\\
y'_k &= V(\hat x_2,\dots, \hat x_k) \hat x_k V^{-1}(\hat x_2,\dots,
\hat x_k) \text{ for } k=3,\dots, n.
\endaligned
$$
It is enough to show that 
$$
z'_k = z_k, \text { and }y'_k=y_k \text{ for } k=2,\dots, n.
$$
For $k=2$ it is obvious.  By Lemma 3.3.1
$$
V(\hat x_2,\dots, \hat x_{k-1}, \hat z)= V(x_1, \dots, x_{k-1},
z)(z-x_1)^{-1},
$$
and by definition $\hat z=(z-x_1) z(z-x_1)^{-1}$.  So, 
$$
z'_k=\{V(x_1,\dots, x_{k-1},z)(z-x_1)^{-1}\}(z-x_1) z(z-x_1)^{-1}\cdot
$$
$$
\cdot\{(z-x_1) V^{-1}(x_1,\dots, x_{k-1}, z)\}= z_k\text{ for }
k=3,\dots, n.
$$
Similarly, $y'_k=y_k \text{ for } k=3,\dots, n\text{ and from } (3.4)$
we have
$$
V(x_1,\dots, x_n, z) = (z_n- y_n)\dots (z_2-y_2)(z_1- y_1)
$$
$\hfill\qed$

\subhead III.1.6 Another expression for coefficients\endsubhead
Another expression for coefficients $a_1,\dots a_n$ in Vieta
decomposition of $V(x_1,\dots, x_n, z)$ may be obtained from Proposition
1.3.7.

\proclaim{Theorem 3.1.4} \text {[GKLLRT]}:
$$
V(x_1,\dots, x_n,z) = z^n+a_1 z^{n-1}+ \dots + a_n,
$$
where for $k=1,\dots, n$
$$
a_k=-\vmatrix
&x_1^n&\dots &x_n^n\\
&x_1^{n-k+1}&\dots& x_n^{n-k+1}\\
&x_1^{n-k-1}&\dots &x_n^{n-k-1}\\
&{}&\dots &{}\\
&1&\dots &1\endvmatrix _{1n}\cdot
\vmatrix
&x_1^{n-1}&\dots&x_n^{n-1}\\
&{}&\dots&{}\\
&x_1^{n-k}&\dots &x_n^{n-k}\\
&{}&\dots&{}\\
&1&\dots&1\endvmatrix^{-1} _{kn}. \tag 3.5
$$
\endproclaim
>From Theorem 3.1.4 we will get Bezout and Vieta formulas expressing
coefficients of the equation 
$$
z^n + a_1 z^{n-1} + \dots + a_n = 0\tag 3.6
$$

{\centerline {\bf \S3.2 Vieta theorem and Bezout theorem}}
\medskip
\subhead III.2.1 Formulas for coefficients\endsubhead  
\proclaim{Lemma 3.2.1} Suppose that $x_1,\dots, x_n$ is an independent
set of roots of the equation (3.6).  Then coefficients $a_1,\dots, a_n$
may be written in the form (3.5)
\endproclaim
\demo{Proof}\enddemo Consider a system of right linear equations 
$$x^n_i + a_1x^{n-1}_i+\dots+ a_{n-1}x_i+ a_n = 0, i=1,\dots, n
$$
with unknowns $a_1,\dots ,a_n$ and use Cramer rules.
$\hfill\qed$

\subhead III.2.2 Bezout and Vieta Theorems\endsubhead
\proclaim{Theorem 3.2.2} \text {(Noncommutative Bezout Theorem)}.  Let
$x_1,\dots,x_n$ be an independent set of roots of equation (3.6).  In
notations of Theorem 3.1.1
$$
z^n+a_1 z^{n-1} + \dots+ a_n=(z_n-y_n)\dots(z-y).
$$
\endproclaim
$\hfill\qed$
\demo{Proof} Use Lemma 3.2.1, Theorem 3.1.4 and Theorem 3.1.1.
\enddemo
\proclaim{Theorem 3.2.3}\text{(Noncommutative Vieta Theorem) [GR3]}
Let $x_1,\dots, x_n$ be an independent set of roots of equation (3.6).
Then coefficients $a_1,\dots, a_n$ of the equation are given by
formulas (3.3).
\endproclaim
\demo{Proof}\enddemo Use Lemma 3.2.1, Theorem 3.1.4 and Theorem 3.1.2
$\hfill\qed$
\bigskip
{\centerline {\bf \S 3.3. Noncommutative symmetric functions}}
\medskip

A general theory of noncommutative symmetric functions was developed
in [GKLLRT]. In fact, in [GKLLRT] were studied different systems of
multiplicative and linear generators in a free algebra ${\Bbb Sym}$
generated by any system of noncommuting variables $\Lambda _i,
i=1,2,\dots $. In [GKLLRT] these variables were called elementary    
symmetric functions but the theory was developed independently
of the origin of $\Lambda _i$'s. 

In this section we apply the general theory to
noncommutative symmetric functions generated by specific $\Lambda _i$'s.
As in the commutative case, they depend of a set of roots of a polynomial
equation.

\subhead III.3.1 A construction of new variables \endsubhead
 We fix $n$ independent indeterminants $x_1,
x_2, \dots, x_n$ and construct new variables $y_1,\dots, y_n$ which are
{\it rational} functions in $x_1, \dots, x_n$:
$$
\aligned
y_1 &= x_1,\\
y_2 &=\vmatrix &x_1 &x_2\\
&1 &1\endvmatrix _{12}x_2 \vmatrix &x_1 &x_2\\
&1 & 1\endvmatrix ^{-1}_{12},\\
&\dots \\
y_n&=\vmatrix &x_1^{n-1} &\dots&x_n^{n-1}\\
&x_1^{n-2}&\dots &x_n^{n-2}\\
&{}&\dots&{}\\ 
&1&\dots&1\endvmatrix _{1n}x_n
\vmatrix &x_1^{n-1}&\dots &x_n^{n-1}\\
&x_1^{n-2} &\dots &x_n^{n-2}\\
&{} &\dots &{}\\
&1{}&\dots&1\endvmatrix ^{-1}_{1n} .
\endaligned$$

In the commutative case $x_i=y_i, i=1,\dots n$. In the noncommutative case
they are obviously different. 

\remark{Remark} Consider a free skew-field $F$ generated by $x_1,\dots , x_n$.
Define on $F$ differential operators $\partial _i$ given by formulas
$\partial _ix_j=\delta _{ij}$ and
satisfying Leibnitz rule $\partial _i (fg)=\partial _i(f)g +
f\partial _i(g)$, for $i=1,\dots n$.

It easy to see that $\partial _iy_j\neq \delta_{ij}$. Consider, however,
$\partial =\partial _1 + \dots +\partial _n$. Then
$$ \partial (V(x_1,\dots , x_k))=0, \ k=2,\dots n$$
and
$$ \partial (y_i)=\partial (x_i)=1, \ i=1,\dots ,n.$$
\endremark
\medskip
\subhead III.3.2. Elementary symmetric functions\endsubhead
\definition{Defintion 3.3.1} The functions
$$
\aligned
\Lambda_1(x_1,\dots, x_n) &= y_1 + y_2\dots + y_n,\\
\Lambda_2(x_1,\dots, x_n) &= \sum_{i<j} y_jy_i,\\
      {}&\dots    {}\\
\Lambda_n(x_1,\dots, x_n) &= y_n\dots y_1\endaligned
$$
are called elementary symmetric functions of $x_1, \dots, x_n.$
\enddefinition
In commutative case these functions are the ordinary elementary
symmetric functions of $x_1,\dots, x_n$.

By Theorem 3.2.3 $\Lambda_i(x_1,\dots, x_n) = (-1)^i a_i, i=1,\dots, n$ where
$x_1, \dots, x_n$ are roots of the equation
$$
x^n+a_1 x^{n-1} + \dots + a_{n-1} x+ a_n = 0.
$$
This implies
\proclaim{Proposition 3.3.2} Functions $\Lambda_i(x_1,\dots, x_n)$ are
symmetric in $x_1, \dots, x_n$.
\endproclaim

\remark{Remark}  The order of $y_1,\dots, y_n$ is essential in the
definition of $\Lambda_i, i=1,\dots, n$.

For example, $\Lambda_2 = y_2y_1$ is symmetric in $x_1, x_2$ but the
product $y_1 y_2$ is not symmetric. 

To see this let $n=2$ and set $d=x_2-x_1$. Then the product
$y_1 y_2$ is symmetric in $x_1, x_2$ if and only if
$x_1d^2=d^2x_1$.
\endremark
\subhead III.2.3 Complete symmetric functions\endsubhead
Let $t$ be a formal variable commuting with $x_i, i=1,\dots n$.
Consider the generating function
$\lambda (t)=1+\Lambda _(x_1,\dots ,x_n)t+\dots 
+\Lambda _n (x_1,\dots,x_n)t^n$.

Following [GKLLRT] define {\it complete} symmetric functions
$S_i(x_1,\dots ,x_n), i=1,2, \dots$ using a generating function

$$1+\sum _i S_i(x_1,\dots , x_n)t^i =\lambda ^{-1}(-t).$$

According to [GKLLRT]
$$
S_k(x_1,\dots, x_n) = \sum_{i_1\leq i_2\leq\dots\leq i_k}y_{i_1}\dots
y_{i_k}, k=1,2,3,\dots  ,\tag 3.7
$$

In commutative case $S_k$'s are just complete symmetric functions and so
we will call $S_k$'s defined by formula (3.7) complete symmetric
functions of $x_1,\dots ,x_n$.

Formula (3.7) shows that $S_k$'s are polynomial in $y_1,\dots ,y_n$
and

\proclaim{Proposition 3.3.3} $S_k(x_1,\dots, x_n)$ are symmetric in
$x_1,\dots, x_n$.
\endproclaim
\remark{Remark} The order of $y's$ in (3.7) is essential:
$
S_2(x_1,x_2) = y^2_1+ y_1y_2 + y_2^2$ is symmetric in $x_1, \dots ,x_n$
and $y_1^2 + y_2y_1 + y_2^2$ is not symmetric. Otherwise, the sum
$y_1^2 + y_2^2$ will be also symmetric and for $n=2$ and
$d=x_2-x_1$ one would have $d^2x_1=x_1d^2$.
\endremark

\subhead III.3.4 Ribbon Schur Functions\endsubhead  
To consider a more general example we need a
vocabulary.  Let $w=a_{i_1}\dots a_{i_k}$ be a word in ordered
letters $a_1<\dots <a_n$.  An integer $m$ is
called a descent of $w$ if $1\leq m\leq k-1$ and $i_m > i_{m+1}$.

Let $J=(j_i,\dots, j_k)$ be a set of positive integers.  
Define the {\it ribbon
Schur function} $R_J(x_1,\dots, x_n)$ by formula
$$
R_J(x_1,\dots, x_n) = \sum y_{i_1}\dots y_{i_m},
$$
where sum is running over all words $w=y_{i_1} \dots y_{i_m}$ whose
descent set is exactly equal to $\{j_1, j_1+j_2,\dots, j_1 + j_2+\dots
+ j_{k-1}\}$.

\proclaim{Theorem 3.3.4} Functions $R_J(x_1,\dots, x_n)$ are symmetric in
$x_1,\dots, x_n$ for any $J$.
\endproclaim
>From [GKLLRT] it follows

\proclaim{Theorem 3.3.5} The set of functions $\{R_J\}$ for all $J$ is
a $\Bbb Q$-linear basis in the free $\Bbb Q$-algebra generated by 
$\Lambda _1(x_1,\dots , x_n), \dots ,\Lambda _n(x_1,\dots ,x_n)$
\endproclaim
\subhead III.3.5 Main Theorem\endsubhead
In a commutative case the main theorem of a theory of commutative
functions says that every symmetric polynomial of $n$ variables
is a polynomial of (elementary) symmetric functions of these
variables. Its analogue for a noncommutative case is given by
\proclaim{Conjecture 3.3.6} Let a polynomial $P(y_1, \dots, y_n)$ over
$\Bbb Q$ be
symmetric as a function of $x_1, \dots, x_n$ then $P(y_1,\dots,
y_n)=Q(\Lambda_1, \dots, \Lambda_n)$ where $Q$ is a noncommutative
polynomial of $n$ variables over $\Bbb Q$.
\endproclaim
We proved this conjecture for $n=2$.

\head IIII Continued fractions and Almost Triangular Matrices\endhead

\subhead IIII.1 Continued fractions and quasideterminants\endsubhead 

Consider an infinite matrix $A$ over a skew-field 
$$
A=\pmatrix &a_{11}&a_{12}&a_{13}&\dots &a_{1n}\dots\\
           &-1    &a_{22}&a_{23}&\dots &a_{2n}\dots\\
           &0     &-1    &a_{33}&\dots &a_{3n}\dots\\
           &0     &0  &-1      &\dots &\dotso\endpmatrix
$$

It was pointed in [GR1], [GR2] that its quasideterminant $|A|_{11}$
can be written as a generalized continued fraction

$$|A|_{11} = a_{11} + \sum_{j_1\neq 1} a_{1j_1}{1\over
a_{2j_1}+\sum \Sb j_2\neq 1,j_1\endSb a_{2j_2} {1\over a_{3j_2}+\dots}}.
$$
Let
$$
A_n=\pmatrix &a_{11} & a_{12} &\dots &a_{1n}\\
&-1 &a_{22} &\dots &a_{2n}\\
&0 &-1 &\dots  &a_{3n}\\
&{} &{} &\dots &{}\\
&\dots  &0 &-1 &a_{nn}\\
\endpmatrix.
$$
The following proposition was formulated in [GR1], [GR2].
\proclaim{Proposition 4.1}
$|A_n|_{11} = P_n Q^{-1}_n$, where
$$
P_n=\sum_{1\leq j_1<\dots < j_k< n} a_{1j_1}a_{j_1+1j_2}
a_{j_2+1j_3}\dots a_{j_k+1n}\tag 4.1
$$
$$
Q_n=\sum_{2\leq j_1<\dots < j_k< n}
a_{2j_1}a_{j_1+1j_2}a_{j_2+1j_3}\dots a_{j_k+1n}.\tag 4.2
$$
\endproclaim
{\bf Proof}. From homological relations one has

$$ 
|A_n|_{11}|A^{1n}_n|^{-1}_{21}= - |A_n|_{1n}|A^{11}_n|^{-1}_{2n}.
$$
We will apply Definition 1.1.3 for a computation of $|A_n|_{1n}$,
$|A^{11}_n|_{2n}$, and $|A^{1n}_n|_{21}$. It is easy to see, that
$|A^{1n}_n|_{21}=-1$. To compute two other quasideterminants we
have to invert triangular matrices. Setting $P_n=|A_n|_{1n}$ and
$Q_n=|A^{11}_n|_{2n}$ we arrive to formulas (4.1), (4.2).
$\square$

\remark{Remark}In a commutative case Proposition 4.1 is well-known.
In this case $P_n=|A_n|_{1n}=(-1)^n\text {det}A_n$
and $Q_n=(-1)^{n-1}\text {det}A^{11}_n$.
\endremark

>From formulas (4.1), (4.2) it follows [GR1], [GR2]
\proclaim {Corollary 4.2}  
The polynomials $P_k$ for $k\geq 0$ and $Q_k$ for $k\geq 1$ are related
via formulas:
$$
P_k=\sum^{k-1}_{s=0} P_s a_{s+1,k}\qquad P_0 = 1,\tag 4.3
$$
$$Q_k=\sum^{k-1}_{s=1} Q_s a_{s+1,k}\qquad Q_1 = 1.\tag 4.4
$$
\endproclaim

The following corollary was pointed out to us by A. Berenstein.
\proclaim{Corollary 4.3}
Consider the matrix $A_n$. Suppose that for any $i\neq j$ and
any $p,q$
$$a_{ij}a_{pq}=a_{pq}a_{ij}$$
and
$$a_{jj}a_{ii}-a_{ii}a_{jj}=a_{ij}\ 1\leq i<j\leq n. $$  
Then
$$ P_n=|A_n|_{1n}=a_{nn}a_{n-1n-1}\dots a_{11} \tag 4.5$$
\endproclaim

The proof follows from formula (4.3). Formulas similar to (4.5)
appeared in [G] and [FK].  

\proclaim{Corollary 4.4} [GR1], [GR2] For a Jacobian matrix
$$
A=\pmatrix &a_1 & 1 & 0 &\dots\\
&-1&a_2 & 1&{}\\
&0 &-1 & a_3 &\dots\endpmatrix
$$
$$|A|_{11} = a_1 + {1\over a_2+{1\over a_3 + \dots}},$$
and
$$
P_k = P_{k-1} a_k + P_{k-2},\ k\geq 2;\ P_0 =1,\ P_1 = a_1,
$$
$$
Q_k = Q_{k+1}a_k + Q_{k-2}.\ k\geq 3; Q_1 = 1,\ Q_2 = a_2.
$$
\endproclaim

In this case $P_k$ is a polynomial of $a_1,\dots ,a_k$ and
$Q_k$ is a polynomial of $a_2,\dots , a_k$.

\subhead IIII.2 Continued fractions and formal series\endsubhead

In notations of the previous subsection infinite continued 
fraction $|A|_{11}$ may be written as a ratio of formal series
of letters $a_{ij}$ and $a^{-1}_{ii}$.

Namely, set
$$
P_\infty=\sum\Sb 1\leq j_1
<j_2\dots<j_k<r-1\\ r=1,2,3,\dots\endSb a_{1j_1}a_{j_1+1j_2}\dots
a_{j_k+1r}a^{-1}_{rr}\cdot\dots\cdot a_{11}^{-1} =
$$

$$
=1+ a_{12}a^{-1}_{22}a^{-1}_{11} +
a_{13}a^{-1}_{33}a^{-1}a^{-1}_{22} a^{-1}_{11} + 
a_{11}a_{23}a^{-1}_{33}a^{-1}_{22}a^{-1}_{11}+\dots ,
$$

$$
Q_\infty=a^{-1}_{11}+\sum\Sb 2\leq j_1
<j_2\dots<j_k<r-1\\ r=2,3\dots\endSb a_{2j_1}a_{j_1+1j_2}\dots
a_{j_k+1r}a^{-1}_{rr}\cdot\dots\cdot a_{11}^{-1} =
$$

$$
=a^{-1}_{11} + a_{23}a^{-1}_{33}a^{-1}_{22}a^{-1}_{11} +
a_{24}a^{-1}_{44}a^{-1}_{33} a^{-1}_{22} a^{-1}_{11} + \dots .
$$
It is easy to see that these formal series are correctly defined.

The following theorem was proved in [PPR]. Here we consider
another proof.
\proclaim{Theorem 4.5}
$$
|A|_{11} = P_\infty\cdot Q^{-1}_\infty.
$$
\endproclaim 

{\bf Proof}. Consider matrix $B=(b_{ij}=a_{ij}a^{-1}_{jj})$,
$i,j=1, 2, 3, \dots $. According to a property of quasideterminants
$|A|_{11}=|B|_{11}a_{11}$.

Applying noncommutative Sylvester Theorem to $B$ with matrix
$(b_{ij}, i,j \geq 3)$ as a pivot one has

$$|B|_{11}= 1+|B^{21}|_{12}|B^{11}|^{-1}_{22}a^{-1}_{11}.$$

It follows
$$
|A|_{11}= (a_{11}|B^{11}|_{22}a^{-1}_{11} 
+ |B^{21}|_{12}a^{-1}_{11})
(|B^{11}|_{22}a^{-1}_{11})^{-1}.\tag 4.6
$$

>From [GKLLRT, Proposition 2.4] it follows that the first
factor in (4.6) equals to $P_{\infty}$, and the second
one to $Q_{\infty}$. $\square $
  
The following application of Theorem 4.5 to Rogers-Ramanujan continued fraction
was given in [PPR].
Consider a continued fraction with two formal variables $x$ and $y$:
$$
A(x,y) = {1\over 1+x{1\over 1+x{1\over 1+\dots}y}y}.
$$
It is easy to see that
$$
A(x,y)=\vmatrix &1 &x &{}&\cdot&{}&{}&{}\\
&-y&1&x&{}&\cdot &{}&O\\
&{}&-y&1&x&{}&\cdot&{}\\
&{}&{}&{}&1&{}&{}&\cdot\\
&{}&O&{}&\ddots&\ddots&{}&\cdot\endvmatrix^{-1}_{11} =
$$

$$
=\vmatrix
&1&x&O&{}&{}\\
&-1&y^{-1}&xy^{-1}&{}&{}\\
&O&-1&y^{-1} &xy^{-1} &{}\\
&{}&{}&-1 &y^{-1}&\ddots\endvmatrix_{11}
$$
Theorem 4.5 implies the following
\proclaim {Corollary 4.6}

$A(x,y) = P\cdot Q^{-1} $, where $Q=yPy^{-1}$ and 
$$
P=1+\sum\Sb k\geq 1\\n_1,\dots,n_k\geq 1\endSb y^{-n_1}x
y^{-n_2}x\dots \, y^{-n_k} x y^{k+n_1+n_2 +\dots+n_k}.
$$
\endproclaim

Following [PPR]  suppose that $xy=qyx$, where $q$ commutes with
$x$ and $y$. Set $z=yx$.

Then Corollary 4.6 implies Rogers-Ramanujan continued fraction identity
$$A(x,y)={1\over
1+{qz\over 1+{q^2 z\over 1+\dots}}}=$$

$$={{1+\sum_{k\geq 1} {q^{k(k+1})\over
(1-q)\dots(1-q^k)}z^k}
\over 1+\sum_{k\geq 1}{q^{k^2}\over
(1-q)\dots(1-q^k)}z^k}.$$
\medskip
\subhead IIII.3 Determinants of almost triangular matrices\endsubhead
\medskip
Consider now a general almost triangualr matrix
$$
B=\pmatrix &a_{11} & a_{12} &\dots &a_{1n}\\
&b_{21} &a_{22} &\dots &a_{2n}\\
&0 &b_{32} &\dots  &a_{3n}\\
&{} &{} &\dots &{}\\
&\dots  &0 &b_{n-1 n} &a_{nn}\\
\endpmatrix.
$$
In this subsection we suppose that all $b_{i+1i}$'s are
invertible and all $a_{ij}$'s are free variables.

Denote by $B_k, k=1,2,\dots ,n-1$ submatrices of $B$ obtained
by deleting first $k$ rows and columns from matrix $B$. 

Set
$$
D(1,2,\dots ,n)=|B|_{11}b^{-1}_{21}|B_1|_{22}b^{-1}_{32}
|B_2|_{33}b^{-1}_{43}\dots b^{-1}_{nn-1}a_{nn}.
 $$
Homological relations imply the following proposition [GR1], [GR2]
(cf. subsection IIII.1).

\proclaim{Proposition 4.7}
$$
(-1)^{n+1}D(1,2,\dots n)=|B|_{1n}=
$$
$$
=\sum_{1\leq j_1<\dots < j_k< n} (-1)^{k+1}a_{1j_1}b^{-1}_{j_1+1j_1}
a_{j_1+1j_2}b^{-1}_{j_2+1j_2}
a_{j_2+1j_3}\dots b^{-1}_{j_k+1j_k}a_{j_k+1n}.
$$
\endproclaim
\remark{Remark}
When $b_{i+1i}=1,\ i=1,\dots ,n-1$ 
and all $a_{ij}$'s commute with each other, then the
 expression $D(1,2,\dots ,n)$
is the determinant of $B$.
\endremark
It is interesting to find out other quasideterminants of an almost
triangular matrix.
\proclaim{Proposition 4.8}\text {[GR1], [GR2]} Set $D(\emptyset )=1$.
 For $i\leq j$
$$
|B|_{ij}=D(1,\dots , i-1)^{-1}D(1,\dots , n)D(j+1,\dots ,n)^{-1}.
$$
\endproclaim


\Refs

\ref\by [C] P.~M. Cohn\book Skew Field Constructions \publ Cambridge
Univ Press \yr 1977 \endref

\ref\by [C1] P.~M. Cohn\book Skew Fields\publ Cambridge Univ Press
\yr 1995 \endref

\ref\by A. Connes and A. Schwarz \paper Matrix Vieta Theorem Revisited
\jour Lett. Math. Physics \yr 1997 \endref

\ref\by [EGR] P. Etingof, I. Gelfand, and V. Retakh\paper Factorization
of Differential Operators, Quasideterminants, and Nonabelian Toda
Field Equations\jour Math. Res. Letters\vol 4\yr 1997\endref

\ref\by [FK] S. Fomin and An. Kirillov\paper Quadratic algebras, Dunkl
operators, and Schubert Calculus \jour preprint\yr 1997 \endref

\ref\by [GR] I. Gelfand and V. Retakh\paper Determinants of Matrices over
Noncommutative Rings \jour Funct. Anal. Appl.\vol 25 \issue 2 \yr 1991
\pages 91-102 \endref

\ref\by [GR1] I. Gelfand and V. Retakh\paper A Theory of Noncommutative
Determinants and Characteristic Functions of Graphs
\jour Funct. Anal. Appl.\vol 26 \issue 4 \yr 1992
\pages 1-20 \endref

\ref\by [GR2] I. Gelfand and V. Retakh\paper A Theory of Noncommutative
Determinants and Characteristic Functions of Graphs. I
\jour Publ. LACIM, UQAM, Montreal \vol 14 \yr 1993 
\pages 1-26 \endref

\ref\by [GR3] I. Gelfand and V. Retakh\paper Noncommutative Vieta
Theorem and Symmetric Functions, in
\book Gelfand Mathematical Seminars 1993-95\publ Birkhauser \yr 1996
\endref

\ref\by [GKLLRT] I. Gelfand, D. Krob, A. Lascoux, B. Leclerc,
V. Retakh, and J.-Y. Thibon\paper Noncommutative Symmetric Functions
\jour Advances in Math\vol 112\issue 2\yr 1995 \pages 218-348 \endref

\ref\by [G] A. Givental\paper Stationary Phase Integrals, Quantum Toda
Lattices, Flag Manifolds and the Mirror Conjecture 
\jour alg-geom 9612001 \yr 1996 \endref

\ref\by [KL] D. Krob and B. Leclerc\paper Minor Identities for
Quasi-Determinants and Quantum Determinants\jour Comm. Math. Phys.
\vol 169 \issue 1\yr 1995 \pages 1-23 \endref

\ref\by [M] A. Molev\paper Noncommutative Symmetric Functions and
Laplace Operators for Classical Lie Algebras\jour Lett. Math. Phys.
\vol 35\issue 2\yr 1995 \pages 135-143 \endref

\ref\by [PPR] I. Pak, A. Postnikov, V. Retakh \paper Noncommutative
Lagrange Inversion \jour preprint \yr 1995 \endref

\ref\by [RS] A. Razumov, M. Saveliev \paper Maximally Nonabelian
Toda Systems \jour Physica D\yr 1997 \endref  
\endRefs
\enddocument